\documentclass[a4paper]{amsart}
\usepackage{amssymb}
\usepackage{color}
\newtheorem{theo}{Theorem}[section]
\newtheorem{prop}{Proposition}[section]

\newtheorem{lemma}{Lemma}[section]
\newtheorem{rema}{Remark}[section]

\numberwithin{equation}{section} 

\def\t{\textrm}

\def\beq{\begin{eqnarray}}
\def\eeq{\end{eqnarray}}
\def\baa{\begin{array}}
\def\eaa{\end{array}}

\newcommand{\bdef}{\begin{definition}}
\newcommand{\be}{\begin{equation}}
\newcommand{\ee}{\end{equation}}
\newcommand{\bt}{\begin{theo}}
\newcommand{\et}{\end{theo}}
\newcommand{\bl}{\begin{lemma}}
\newcommand{\el}{\end{lemma}}
\newcommand{\bp}{\begin{prop}}
\newcommand{\epr}{\end{prop}}

\def\dsp{\displaystyle}
\def\pa{\partial}
\def\ep{\epsilon}
\def\ov{\overline}

\let\dsp=\displaystyle
\def\R{{\mathbb R}}

\def\varpsi{\psi}
\def\li{\lambda_i}
\def\lai{\lambda_i}
\def\v{v_i}
\def\u{u_i}
\def\vx{{v_i}_x}
\def\ux{{u_i}_x}
\def\uxx{{u_i}_{xx}}
\def\vt{{v_i}_t}
\def\ut{{u_i}_t}

\def\beti{\beta_i}
\def\Di{D_i}
\def\phii{\phi_i}
\def\phx{{\phi_i}_x}
\def\pht{{\phi_i}_t}
\def\phxx{{\phi_i}_{xx}}
\def\phtx{{\phi_i}_{tx}}

\def\pho{\phi_{0i}}

\def\uo{u_{0i}}
\def\vo{v_{0i}}

\def\vox{{v_{0i}}_x}
\def\M{{\mathcal{M}}}
\def\U{\mathcal{O}}
\def\E{\mathcal{I}}
\def\A{\mathcal{A}}
\def\Ii{I_i}

\def\iIi{\int_{\Ii}}
\def\nld{\Vert_{2}}

\def\nhu{\Vert_{H^1}}

\def\soT{\sup_{[0,T]}}
\def\iot{\int_0^t}

\def\nhd{\Vert_{H^2}}
\def\ioT{\int_0^T}

\def\suE{\sum_{i\in \E}}
\def\suM{\sum_{i\in \M}}
\def\suU{\sum_{i\in \U}}

\def\suEn{\sum_{i\in \E^\nu}}
\def\suUn{\sum_{i\in \U^\nu}}

\def\iota{\int_0^\tau}
\def\ioTa{\int_0^\tau}
\def\del{\Delta^h}
\def\T{\mathcal{T}}

\begin{document}
\title[A hyperbolic chemotaxis model on a network]
{Global smooth solutions for a hyperbolic chemotaxis  model on a network}
\author[ F. R. Guarguaglini and R. Natalini]{ F. R. 
Guarguaglini$^\dag$  and
R. Natalini$^\diamond$}
 
\thanks{\noindent 
$\dag$ Dipartimento di Ingegneria e Scienze dell'Informazione e Matematica, Universit\`a degli
Studi di L'Aquila, Via Vetoio, I--67100 Coppito (L'Aquila), Italy. E--mail:
guarguag@univaq.it \\
\indent $\diamond$ Istituto per le Applicazioni del Calcolo ``M.
Picone'',
Consiglio Nazionale delle Ricerche, 
Rome, Italy. E--mail: roberto.natalini@cnr.it}
\date{}
\begin{abstract} 
In this paper we study a semilinear hyperbolic-parabolic system modeling  biological
phenomena evolving on a network composed by  oriented arcs. We prove the existence of global (in time) smooth solutions to this problem. The result  is obtained by using energy estimates with suitable transmission conditions at nodes.
\end{abstract}
\maketitle
\bigskip

\begin{section}
{\bf Introduction}

In this paper we consider a semilinear hyperbolic-parabolic system which describes  chemosensitive  movements of bacteria, cells or other microorganisms on an artificial scaffold. 

A mathematical approach to   phenomena involving chemotaxis
has been largely developed in the last thirty years, see \cite{horst,perth,mur}, mostly by means of the 
study of the Patlak-Keller-Segel system \cite{KS}; this model was constituted by a parabolic  equation governing the evolution of the density of  cells, and a parabolic or elliptic one for the evolution of concentration of chemoattractant. 

On the contrary, here, following \cite{AG,gumanari}, we consider a model where the evolution of the density of  cells
is described by a hyperbolic system, coupled with a parabolic equation for the chemoattractant, which 
in one space dimension reads
\be\label{eq}\left\{
\begin{array}{ll}
u_t+\lambda v_x=0\ ,
\\ \\
v_t +\lambda u_x=\phi_x u-\beta v\ ,\qquad\qquad x\in I\, ,t\geq 0\ ,
\\ \\
\phi_t=D\phi_{xx}+au-b\phi\ \ .
\end{array}\right.\ee
The unknown $u$ stands for  the concentration of  cells, $v$ is their average flux and $\phi$ is the
concentration of chemoattractant produced by the cells themselves; the source term $\phi_x u$ is the nonlinear chemotactic term. As regard to the various parameters, $D>0$ is the diffusion coefficient of chemoattractant, $a\geq 0$ and $b>0$ are respectively its production and degradation rates, $\beta>0$ is the friction coefficient of substrate; finally,  each individual can move at a constant velocity, whose modulus is $\lambda\geq 0$, towards right or left along the axis.

Hyperbolic models have been recently introduced \cite{DH,FLP,Hillen,perth}   since they yield a more realistic finite speed of propagation, in contrast with the parabolic ones, and allow a better observation of the phenomena during the transitory. Models like (\ref{eq}) were proposed  in \cite{Segel,AG}, introducing the chemotactic term in the Cattaneo equation,
and later their solutions were studied in \cite{hillerho,hillenst,gumanari}.

Here we consider the one dimensional model (\ref{eq}) cast on a network formed by $n$ nodes $N_\nu$ and $m$ oriented arcs connecting the nodes, $I_i$; moreover, each arc is characterized by a typical velocity $\lambda_i$ and a diffusion coefficient $D_i$ . 
On each arc $I_i$ we consider the triple of unknowns $(\u,\v,\phii)$.

This model arises as a preliminary tool in tissue-engineering research, to describe the process of dermal wound healing: the fibroblasts are seeded on a polymeric 
scaffold and they move along it to fill the wound, driven by chemotaxis \cite{harley,mandal,spadac}. In our mathematical  model the arcs mimic the fibers of the scaffold and they intersect at nodes and $\u,\phi_i$
are the densities of fibroblasts and  chemoattractant on each of them.

The aim of this paper is to give a proof of global existence of smooth solutions to problem (\ref{eq}) complemented with initial, boundary and transmission conditions at nodes.
The rigorous statement of the problem is presented in Section 2.

Since the Cauchy and the Neumann problems were considered in \cite {gumanari}, the crucial point of our study is in  the transmission conditions, which heavily characterize  the problem, since they are the coupling among the solutions on different arcs.
 So far,  our analysis is based on the formulation of suitable transmission conditions at nodes, also in the   preliminary proof of local existence. More precisely, we impose transmission conditions which guarantee the conservation of the fluxes, both for the hyperbolic system and the parabolic equation; in this way  the energy of the linearized homogeneous version of the system decades in time.
Section 3 is devoted the motivation and the derivation of these conditions.

We mention that hyperbolic models on networks  have been previously studied in \cite{pic,zua1,zua2}; moreover a parabolic chemotaxis model on networks was studied in \cite {4}, with continuity conditions at node. Finally, in
 \cite{BNR} the authors treat the same our model from a numerical point of view, with slightly different transmission conditions.

In Section 4 we prove the first result of this paper, namely the existence and uniqueness of local solutions. This result is obtained by means of the linear contraction semigroups theory coupled with the abstract theory of nonhomogeneous and semilinear evolution problems ; in fact the right position of transmission
conditions  allows to use the properties of m-dissipative linear operators.

Finally, in Section 5, we prove the existence of global solutions in the case of small (in a suitable norm) initial data . The proof of this result needs some further conditions at node
providing the tools for controlling the growth of the unknowns' traces.

\end{section}
\begin{section} 
{\bf Statement of the problem}

We consider a connected graph $G=(\mathcal N, \mathcal A)$ composed by a 
set  $\mathcal N$ of $n$ nodes, $\mathcal N=\{ N_\nu:\nu\in \mathcal P=\{1,2,...,n\}\}$
 and a set $\mathcal A$
of $m$ oriented arcs, $\mathcal A=\{I_i:i\in \mathcal M=\{1,2,...,m\}\}$. Moreover we denote by $a_j$ , $j\in \mathcal J=\{1,2,...,l\}$  the external points of the graph. 

Since we are on an oriented network, for each node $N_\nu$ we consider the set of incoming arcs 
$\mathcal A_{in}^\nu=\{I_i:i\in \mathcal I^\nu\}$
and the set of the outgoing ones
$\mathcal A_{out}^\nu=\{I_i:i\in \mathcal O^\nu\}$; let $\mathcal M^\nu=\mathcal I^\nu\cup\mathcal O^\nu$.

Each oriented arc $I_i$ is a compact one dimensional interval; if $I_i$ is an  external arc, connecting a boundary point
$a_j$ to a node $N_\nu$, then it has the form $I_i=[a_j,N_\nu]$ when  $i\in\mathcal I^\nu$ and the the form $I_i=[N_\nu, a_j]$ when  $i\in\mathcal O^\nu$; if $I_i$ is an internal arc, connecting the nodes $N_\nu$ and $N_\mu$ then it has form $I_i=[N_\nu,N_\mu]$ if $i\in\mathcal O^\nu\cap\mathcal I^\mu$.

A function $f$ defined on $\mathcal A$ is a m-pla of functions $f_i$, $i\in \mathcal M$, each one defined 
on $I_i$; $\dsp L^2(\mathcal A)=\cup_{i\in\M}L^2(I_i)$, 
$\dsp H^s(\mathcal A)=\cup_{i\in\M}H^s(I_i)$
 and
$$\Vert f\Vert_2:=\suM \Vert f_i\Vert_2\ , \ 
\Vert f\Vert_{H^s}:=\suM \Vert f_i\Vert_{H^s}.$$ 

We consider the evolution on the graph $G$ of the following one dimensional problem 

\begin{equation}\label{sysi}
    \left\{ \begin{array}{l}
\partial_t \u +\lambda_i \partial_x \v=0\ ,
\\ \\
\partial_t \v +\lambda_i \partial_x \u= \u \partial_x \phi_i -\beta_i \v\ ,
\qquad\qquad x\in I_i\, ,t\geq 0\, , i\in \M,
\\ \\
\partial_t \phii= D_i \partial_{xx} \phi_i +a\u -b\phi_i\ ,
\end{array}\right.\end{equation}
where $\lambda_i, a\geq 0\ ,\ b,D_i,\beta, >0 .$  Actually, the coefficients $a$ and $b$ could depend on the arc $I_i$ in consideration, but $\dsp\frac {a_i}{b_i}$ should be  constant for $i\in\M$.

We complement the system with the initial conditions
\be\label{uvic}
u_{i0},v_{i0}\in H^1(I_i)\ ,\ 
\phi_{i0}\in H^2(I_i)
\ \t{ for } i\in\M\ .\ee

As regard to the boundary conditions, at each outer point $a_i$ we set the null flux conditions

\be\label{vbc}
\v(a_i, t)=0\ ,\ 
\qquad 
\  t>0\ ,
\ee
\be\label{phibc}
\phx(a_i, t)=0\qquad 
\ t>0\ .\ee



Besides, at each node $N_\nu$ we impose the following transmission conditions for $\phi_i(N_\nu,t)$
\be\label{KC}\left\{\begin{array}{ll}\displaystyle
\Di \phx(N_\nu,t)=  \sum_{j\in\M^\nu} \alpha_{ij}^\nu(\phi_j(N_\nu,t)-\phii(N_\nu,t))\  \ ,\ i\in \E^\nu \ ,\  t> 0\ ,
\\ \\ \dsp
-\Di \phx(N_\nu,t)=  \sum_{j\in\M^\nu} \alpha_{ij}^\nu(\phi_j(N_\nu,t)-\phii(N_\nu,t))\  ,\ \  i\in \U^\nu \ ,\  t> 0\ ,
\\ \\
\alpha_{ij}^\nu\geq 0\ ,\ \alpha_{ij}^\nu=\alpha_{ji}^\nu \ \t{ for all } i,j\in \M^\nu\ ,
\end{array}\right.\ee
which imply the continuity of the flux  at  node, for all $t>0$,
$$\suEn D_i\phx(N_\nu,t)=\suUn D_i \phx(N_\nu,t)\ .$$

In similar way we impose some transmission conditions for the unknowns $\v(N_\nu,t)$ and $\u(N_\nu,t)$
\be\label{TC}\left\{\begin{array}{ll}\dsp
-\lai\v(N_\nu,t)=\sum_{j\in\M^\nu}
K_{ij}^\nu \left(u_j(N_\nu,t)-\u(N_\nu,t)\right)\ ,\ i\in \E^\nu\ ,\ t>0\ ,\\ \\ \dsp
\lai\v(N_\nu,t)=\sum_{j\in\M^\nu}
K_{ij}^\nu \left(u_j(N_\nu,t)-\u(N_\nu,t)\right)\ ,\  i\in\U^\nu \ , \ t>0\ ,\\ \\
K_{ij}^\nu\geq 0\ ,\ K_{ij}^\nu=K_{ji}^\nu \ \t{ for all } i,j\in \M^\nu\ .
\end{array}\right.\ 
\ee
The above conditions ensure the conservation of the 
flux of the density of cells at each node $N_\nu$ , for $t>0$,

$$
\suEn \lai\v(N_\nu,t)=\suUn\lai\v(N_\nu,t)\ .
$$

We notice that the previous equality corresponds to the conservation of the total mass
$$\suM\iIi\u(x,t)\ dx=\suM\iIi\uo(x)\ dx\ ,$$
which means that no death nor birth of individuals occurs during the observation.

The constrains on the coefficients in the transmission conditions will be widely motivated in the next section.

Finally, we impose the following compatibility conditions
 \be\label{compat} 
u_{i0},v_{i0}, \phi_{i0}\t{  satisfy conditions (\ref{vbc})-(\ref{TC})  
for all } i\in\M\ .
\ee

Provided the above conditions, we will prove the existence and uniqueness of a local solution
 to problem (\ref{sysi})-(\ref{compat}),
$$(u,v)\in C([0,T];H^1(\mathcal A))\cap C^1([0,T],L^2(\mathcal A)\ ,\ 
\phi \in C([0,T];H^2(\mathcal A))\cap C^1([0,T],L^2(\mathcal A)\ .$$

The final aim of the paper is the study of existence of global solutions under the assumption
of smallness of the data. This question involves the proof of estimates for the traces of the unknowns $\u$ at nodes $N_\nu$, which can be derived under the  further assumption that, for all $\nu\in\mathcal P$,  for some $k\in\M^\nu$, 
the coefficients $K_{ik}^\nu$ are non null for all $i\in\M^\nu$, $i\neq k$.

\end{section}

\bigskip

\begin{section}
{\bf Transmission conditions}

In the present section we explain the derivation of the transmission conditions (\ref{KC}) and (\ref{TC})
at nodes. Such conditions have to guarantee two properties of the model; first, the conservation of the flux of the  density of cells, at each node $N_\nu$ 
\be\label{FC}
\suEn\lambda_i\v(N_\nu,t)=\suUn\lambda_i\v(N_\nu,t)
\ee 
and   the conservation of the flux of the parabolic equation
\be\label{FC2}\suEn D_i\phx(N_\nu,t)=\suUn D_i \phx(N_\nu,t)
\ ,\ee
for $t>0$.
Moreover, it is necessary to deal with  dissipative conditions; in other words, the sum of the $m$ energies of the linear version of the hyperbolic systems,  
\begin{equation}\label{sysi2}
    \left\{ \begin{array}{l}
\partial_t \u +\lambda_i \partial_x \v=0
\\ \\
\partial_t \v +\lambda_i \partial_x \u= -\beta_i \v
\end{array}\right.\ee
and the sum of the ones  of the homogeneous parabolic equations 
\be\label{fi2}\partial_t \phii= D_i \partial_{xx} \phi_i -b\phii\ee
have to decay in time.  It is clear, by using the integration by parts,  that the dissipation of such energies
$$
E_1(t)=\suM\iIi \left(\u^2(x,t)+\v^2(x,t)\right)\, dx\ ,\ 
E_2(t)=\suM\iIi \phii^2(x,t) \, dx\, ,
$$
is strictly linked with the following sign properties  of the terms at nodes
\be\label{n1} \Gamma_1^\nu(t)=\suEn \lambda_i\v \u (N_\nu,t)-\suUn\lambda_i\v\u(N_\nu,t)\geq 0\ ,\  \nu\in\mathcal P\ ,
\ee
\be\label{n2} \Gamma_2^\nu(t)=\suEn D_i \phii\phx(N_\nu,t)-\suUn D_i\phii\phx(N_\nu,t)\geq 0\ ,\ \nu\in\mathcal P\ .\ee
In particular, the above conditions ensure that the linear unbounded operators appearing in
the problems (\ref{sysi2}), (\ref{fi2}) are dissipative; this property is crucial  to apply the theory of linear contraction semigroups,  to prove the existence of local solutions, in the next section.

In order to derive transmission conditions which imply 
(\ref{n1}), (\ref{n2}), for some $\nu\in\mathcal P$ , we
first  consider the simple case of two arcs, $I_1$ incoming in $N_\nu$ and $I_2$ outgoing from $N_\nu$.
Here the flux conservation condition (\ref{FC}), together
with inequality (\ref{n1}), reads
$$\lambda_1 v_1(N_\nu,t)(u_1-u_2)(N_\nu,t)= \lambda_2 v_2(N_\nu,t)(u_1-u_2)(N_\nu,t)\geq 0\ ;$$
a possible condition to make this inequality true, is to link the values $\v(N_\nu,t)$ and $\u(N_\nu,t)$ 
as follows:
$$\lambda_1 v_1(N_\nu,t)=\lambda_2 v_2(N_\nu,t)=K^\nu(u_1(N_\nu,t)-u_2(N_\nu,t))\ ,\ \t{ for some } K^\nu\geq 0
\ .$$

Arguing in similar way, using (\ref{FC2}) and (\ref{n2}),  we obtain the transmission conditions for $\phi$ at node
$$D_1 \phi_{x1}(N_\nu,t)=D_2 \phi_{x2}(N_\nu,t)=\alpha^\nu (\phi_2-\phi_1)(N_\nu,t) \ ,\ \t{ for some } \alpha^\nu\geq 0\ .$$

Let us also notice that such kind of conditions for parabolic equations, were introduced in \cite{KK}  in the description of passive transport through biological membranes and they 
are known as Kedem- Katchalsky permeability conditions.

In the case of $m$ arcs intersecting in $N_\nu$, 
the continuity of the flux (\ref{FC}) and the inequality (\ref{n1}) provide the following conditions at node, for $j\neq i$
$$\dsp-\suEn \lambda_i\v(N_\nu,t)(u_j-\u)(N_\nu,t)+\suUn\lambda_i\v(N_\nu,t)(u_j-\u)(N_\nu,t)\geq 0
\  \t{ for all } j\in \M^\nu ;$$
hence some  relations among the values $\v(N_\nu,t)$ and the jumps $(u_j-\u)(N_\nu,t)$, for $i,j\in\M^\nu$, are espected to be asked.

We assume that the terms $v_i(N_\nu,t)$ are
linear combinations of the jumps $\u(N_\nu,t)-u_j(N_\nu,t)$, $j\in M^\nu$:
\begin{equation}\label{vN}\begin{array}{ll}\dsp
-\lai \v(N_\nu,t) =\sum_{j\in\M^\nu}K_{ij}^\nu (u_j(N_\nu,t)-\u(N_\nu,t))\ ,\qquad i\in \E^\nu\ ,\\ \\ \dsp
\lai \v(N_\nu,t) =\sum_{j\in\M^\nu} K_{ij}^\nu (u_j(N_\nu,t)-\u(N_\nu,t))\ ,\qquad i\in \U^\nu\ .\end{array}\ee

Inserting the above positions in the flux continuity relation we obtain

$$\dsp\sum_{i,j\in\M^\nu} K_{ij}^\nu \left(\u -u_j\right)(N_\nu,t) =0\ ;$$
so, we obtain a first constrain to the coefficients
$$\dsp\sum_{i\in\M^\nu} \left(K_{ij}^\nu-K_{ji}^\nu\right) =0\quad \t{ for all } j\in \M^\nu\ .$$

Now we consider the dissipation condition (\ref{n1}) which reads
$$\dsp\dsp\sum_{i,j\in\M^\nu} K_{ij}^\nu \left(\u -u_j\right)(N_\nu,t) \u(N_\nu,t)\geq 0\ ;$$
 sufficient conditions to guarantee the above inequality are
$$K_{ij}^\nu=K_{ji}^\nu\, , K_{ij}^\nu\geq 0 \t{ for all } i,j\in \M^\nu\ .$$

The corresponding conditions on the coefficients $\alpha_{ij}^\nu$ follow by similar calculations.

Finally, we notice that in \cite{BNR} the authors study our problem by a numerical point of  view, introducing transmission conditions for the Riemann invariants of the hyperbolic part of the system, $w^\pm_i=\frac 1 2 (u_i\pm v_i)$, which, in some cases, are equivalent to the present ones.

\end{section}

\begin{section}
{\bf Local  existence }

Let $\dsp Y=\cup_{i\in\M} (L^2(I_i))^2$ $X=L^2(\mathcal A)$.

We consider the linear operator  $A_1:D(A_1)\to Y$, 
\be \begin{array}{ll}
D(A_1)= \{U=(u,v)\in \cup_{i\in\M} (H^1(I_i))^2
: (\ref{vbc}),(\ref{TC}) \t{ hold }\}\\ \\
 A_1 U=\left\{(-\lai \vx, -\lai \ux)\right\}_{i\in\M}\ ,
\end{array}\ee
and the linear operator
$A_2:D(A_2)\to X$,
\be \begin{array}{ll}
D(A_2)= \left\{\phi\in H^2 (\mathcal A): (\ref{phibc}), (\ref{KC}) \t{ hold }\right\}\\ \\
A_2(\phi)=\left\{D_i\phxx-b\phii\right\}_{i\in\M}\ .

\end{array}\ee

We will obtain the existence of local solutions to problem (\ref{sysi})-(\ref{compat}) by  the fixed point technique, combining the local solutions of the two disjointed problems

\be\label{A1} \left\{\begin{array}{ll}
U\in C([0,T];D(A_1))\cap C^1([0,T];Y) \\ \\
U'(t) =A_1 U(t) +F(t,U(t))\ ,\quad t\in[0,T]\ ,\\ \\
U(0)=(u_0,v_0)\in D(A_1)\ ,
\eaa\right.\ee
where
$F(t,U(t))=\{(0,f_i(t)\, u_i(t)-\beta_i v_i(t))\}_{i\in \M}$ and $f$ is a suitable given function  to be specified below, and

\be\label{A2} \left\{\begin{array}{ll}
 \phi \in C([0,T];D(A_2))\cap C^1([0,T];X) \\ \\
\phi'(t) =A_2 \phi(t) + g(t))\ ,\quad t\in[0,T]\ ,\\ \\
\phi(0)=\phi_0\in D(A_2)\ ,
\eaa\right.\ee
where $g$ is a suitable given function  to be specified below.

We have some results concerning the solutions of such problems. In order to simplify the notations,  we will give the proofs in the  case
of a graph composed by a single node $N$ and $m$ arcs $I_i$ connecting that node to the external points $a_i$, 
$i\in\M=\{1,2,...,m\}$. In the general case there are no major differences when treating integrals on the external arcs; on the other hand,  two transmission terms arise when integrating on internal arcs $I_i=(N_\nu,N_\mu)$, each one corresponding to a  node. Then, the sum of all the transmission terms at each node of the graph  can be treated  separately, as in the case of a single node.

\bp\label{hA} Let $T<1$, 
let $g\in C([0,T];H^1(\mathcal A))\cap C^1([0,T],L^2(\mathcal A))$ , 
$\dsp M>\soT\Vert g(t)\nhu$ and  $K> \Vert \phi_0\nhd+4 M$;
then there  exists a unique solution 
to problem (\ref{A2})
and
$$\sup_{t\in[0,T]}\Vert \phi(t)\nhd \leq K\ .$$
Moreover, $\phi\in H^1((0,T);H^1(\mathcal A))$. 
\epr

\begin{proof}
We are going to prove that $A_2$ generates a semigroup in $X$.

First, $A_2$ is dissipative in $X$:
\be\label{diss}\begin{array}{ll}\dsp
(A_2\phi,\phi)= \suM \iIi \left(D_i\phii\pa_{xx}\phii  -b\phii^2\right) \ dx
=-\suM
\iIi \left(D_i \phx^2 + b \phii^2\right) \ dx
\\ \\ \dsp
+ \suE D_i \phx(N,t)\phii(N,t) -\suU D_i\phx(N,t)\phii(N,t) 
\\ \\ \dsp
=-\frac 1 2\sum_{i,j\in M} \alpha_{ij} (\phi_j(N,t)-\phii(N,t))^2 -
\iIi \left(D_i \phx^2 + b \phii^2\right)\ dx\ ,
\end{array}\ee
where we used the trasmission conditions (\ref{KC}).

Moreover, for all  $\varphi\in L^2(\mathcal A)$  there exists  $\phi\in D(A_2)$  such that  $\phi-A_2\phi= \varphi$, i.e.
$A_2$ is a m-dissipative operator in $X$.
In order to prove this fact we introduce the bilinear form $a(\phi,\psi):(H^1(\mathcal A))^2 \to \R$
\be \begin{array} {ll}\dsp
a(\phi,\psi)=\suM\iIi \left(D_i\phx\psi_{ix} +(1+b) \phii\psi_i\right)\ dx
\\ \\ \dsp
-\sum_{i,j\in M}\alpha_{ij}\left(\phi_j(N) -\phii(N)\right) \psi_i(N)\ ;
\end{array}\ee
it is easy to verify that the form is continuous and coercive. Then,
by the Lax-Milgram theorem, we know  that, for all $\varphi\in L^2(\mathcal A)$, there exists a unique
$\phi\in H^1(\mathcal A)$ such that, for all $\psi\in H^1(\mathcal A)$ it holds
$$a(\phi,\psi)=\suM\iIi \varphi_i\psi_i\ dx ;$$
taking $\psi_i\in H^1_0(I_i)$ for all $i\in \M$, we obtain that 
$\phx\in H^1(I_i)$, then
\be\begin{array}{ll}\dsp
\suM\iIi (-D_i\phxx +(1+b) \phii)\psi_i \ dx 
+\suE D_i\left(\phx(N)\psi_i(N) -\phx(a_i)\psi_i(a_i)\right) \\ \\
\dsp
-\suU D_i\left(\phx(N)\psi_i(N) -\phx(a_i)\psi_i(a_i)\right)
 -\sum_{ij\in M} \alpha_{ij}(\phi_j-\phii)(N)) \psi_i(N)\\ \\
\dsp =\suM\iIi \varphi_i\psi_i\ dx\ .\end{array}\ee
The above equality holds for all $\psi_i\in C^\infty_0(I_i)$, then 
$$-\phxx +(1+b) \phii =\varphi_i\quad a.e.\t{ for all }i\in \M$$
and moreover, thanks to suitable choices of $\psi_i(N), \psi_i(a_i)$,
we obtain that $\phi$ satisfies the right boundary conditions to belong to $D(A_2)$.

Then $A_2$ is m-dissipative \cite{CH} and generates a contraction semigroup $\T_2(t)$ in 
$X$; since $g\in C^1([0,T],L^2(\mathcal A))$ we can apply the theory for nonhomogeneous problems in \cite{CH}  to conclude that  there exists a unique  solution to the  problem (\ref{A2}),
given by
$$\phi(t) =\T_2(t)\phi_0 +\iot \T_2(t-s) g(s)\ ds \ .$$
We set 
$$\dsp\mathcal F(t):=\iot \T_2 (t-s) g(s)\ ds\ .$$ $\mathcal F \in C^1([0,T]; L^2(\mathcal A))$ and
\be\label{Fpr}  \mathcal F'(t)=\iot \T_2(s) g'(t-s) \ ds +\T_2(t) g(0) \ ;\ee
moreover $\mathcal F \in C([0,T]; D(A_2))$ and $A_2 \mathcal F(t)=\mathcal F'(t)-g(t)$, see
\cite{CH}.

Hence
\be\begin{array}{ll}\dsp
\Vert \phi(t)\Vert_{D(A_2)} \leq \Vert \phi_0\Vert_{D(A_2)} +
\Vert \mathcal F(t) \Vert_{X}+\Vert A_2\mathcal  F(t) \Vert_{X}\\ \\ \dsp
\leq \Vert \phi_0\Vert_{D(A_2)} +\iot \Vert g(s)\Vert_X\ ds +\Vert \mathcal F'(t)\Vert_X
+\Vert g(t)\Vert_X .
\end{array}\ee
Now, using (\ref{Fpr}) we have
\be\label{contr}\begin{array}{ll}
 \dsp\Vert \phi(t)\Vert_{D(A_2)} \leq \Vert \phi_0\Vert_{D(A_2)} +\Vert g(0)\Vert_X+
\Vert g(t)\Vert_X \\ \\ \dsp+ t\left(\sup_{s\in[0,t]} \Vert g'(s)\Vert_X +
\sup_{s\in[0,t]} \Vert g(s)\Vert_X \right)\end{array}\ee
and, thanks to the condition on $T$, we have
$$\dsp\sup_{t\in[0,T]} \Vert \phi(t)\Vert_{H^2}\leq K\ ,$$
where $K$ is the quantity introduced in the statement of the theorem.

To prove the last claim, it is sufficient to prove that there exists  $C>0$ such that, for all $0<t_1<t_2<T$,
\be\label{w12}\int_{t_1}^{t_2}\Vert \phi_x(t+h)-\phi_x(t)\Vert_{2}^2\, dt \leq C |h|^2\ ,\ee
for all $h\in\R$ with $\vert h\vert <\min\{t_1,T-t_2\}$.

Let $\Delta^h \psi(t):= \psi(t+h)-\psi(t)$; using the equation we can write
$$\int_{t_1}^{t_2}\iIi \left((\Delta^h \phii)_t \Delta^h \phii - D_i (\Delta^h \phii)_{xx}\Delta^h\phii +
\Delta^h g_i\Delta^h \phii-(\Delta^h \phii)^2 \right) \, dx dt =0\ ;$$
then we have
\be\begin{array}{ll}\dsp
\suM \left(\iIi(\Delta^h \phii(t_2))^2 \, dx+
\int_{t_1}^{t_2}\iIi (\Delta^h\phx)^2\ dx dt\right) \\ \\ \leq  \dsp
C \int_{t_1}^{t_2} \left(\suE D_i(\Delta^h \phx) (\Delta^h \phii) (N,t) -
\suU D_i (\Delta^h \phx) (\Delta^h \phii) (N,t)\right) dt\\ \\
\dsp + C\suM\left(\iIi (\Delta^h \phii(t_1))^2 dx+
\int_{t_1}^{t_2}\iIi (\Delta^h g_i)^2 dx dt\right)\ ,
\end{array}\ee
hence inequality (\ref{w12}) follows thanks to nonpositivity of the first term on the
right hand side (as in (\ref{diss})), since $\phi,g\in C^1((0,T);L^2(\mathcal A))$.

\end{proof}

In order to treat problem (\ref{A1}) and to prove the results of the next proposition, we need
the following lemma.

\bl\label{lemma}
Let $ W=(w,z)\in \dsp\cup_{i\in\M} (C^\infty_0(I_i))^2$; there exists a unique $U=(u,v)\in D(A_1)$ such that
$(I-A_1)U=W$.
\el
\begin{proof}
Let $\theta_i=1$ if $i\in\E$ and $\theta_i=-1$ if $i\in\U$. 
We consider the elliptic problem
\be
\left\{
\begin{array}{ll}\dsp - \li^2 \uxx +\u =-\li z_{ix}+ w_i\\ \\ \dsp
\ux(a_i) =0\ , \quad \theta_i\li^2\ux(N)=\sum_{j\in \M} K_{ij} (u_j(N)-\u(N))\ ;\end{array}
\right.
\ee
in the proof of the previous proposition, in the steps to obtain the $m-$ dissipativity of $A_2$, we showed that such problem has a unique solution $u$ whose components $\u$, in the present case, belong to $C^\infty(I_i)$.
Now we set $\v=-z_i-\li \ux$; then $\v\in C^\infty(I_i)$, $\li \vx +\u=w_i$ and $\v(a_i)=0$, 
$\dsp -\theta_i\lambda_i\v(N)=\sum_{j\in M} K_{ij} (u_j(N)-\u(N))$.

\end{proof}

Notice that, if $f\in C([0,T];H^1(\mathcal A))$,   then $F(t,U(t))=f(t)u(t)-\beta v(t)$ is a globally Lipschitz function in
$E:= \dsp \cup_{i\in \M} (H^1(I_i))^2$, with Lipschitz costant $\dsp L_F=L_F\left(\sup_{t\in [0,T]}\Vert f(t)\nhu\right)$; more precisely
$$\soT\Vert F(t,U_1(t))-F(t,U_2(t) \Vert_E \leq L_F 
\soT\Vert U_1(t)-U_2(t) \Vert_E\ .$$

\bp\label{41}
Let $f\in C([0,T_1];H^1(\mathcal A))\cap H^1((0,T_1);L^2(\mathcal A))$,
 $K>\dsp\sup_{[0,T_1]} \Vert f(t)\nhu$,
$\dsp M>2(\Vert u_0\nhu +\Vert v_0\nhu)$ and 
$T<\min\{T_1,(2 L_F(K))^{-1}\}$;
then there  exists a unique solution 
to problem (\ref{A1}) and
$$\sup_{t\in[0,T]}\Vert U(t)
\Vert_E
\leq M\ .$$
\epr

\begin{proof}
First we prove that $A_1$ generates a contraction semigroup in $Y$. 
$A_1$ is a dissipative operator in $Y$: let $U\in D(A_1)$ 
\be\begin{array}{ll}\dsp
(A_1 U, U)=\suM \iIi(-\lai \vx \u-\lai \ux \v) \\ \\
=-\dsp
\left[\suE\lai\v(N) \u(N) -\suU\lai\v(N) \u(N)\right]  \ ;\end{array}
\ee
now,  using the transmission conditions (\ref{TC}) we have
\be\begin{array}{ll}\dsp
(A_1 U, U)=-\frac 1 2\sum_{i,j\in\M} K_{ij}
(u_j(N)-\u(N))^2\ .\end{array}
\ee 
 In order to prove that $A_1$ is a m-dissipative operator in $Y$, we introduce the bilinear form
$a:D(A_1)\times D(A_1)\to \R$
$$a(U,\ov U)=\suM\iIi \left( (\li \vx+\u)(\li \ov\vx+\ov\u) + (\li \ux+\v)(\li \ov\ux+\ov\v) \right)\,dx\ ;$$
$a$ is continuous and , using the boundary and transmission conditions, it is coercive:
$$a(U,U)=\suM\iIi\left( \li^2 \vx^2 +\u^2+\li^2 \ux^2+v^2 \right)\, dx+ \sum_{i,j\in\M} K_{ij}
(u_j(N)-\u(N))^2\ \ .$$
Thanks to the Lax-Milgram theorem, for all $\Psi=(\varpsi_1,\varpsi_2)\in (L^2(\mathcal A))^2$,  there exists
a unique $U\in D(A_1)$ such that, for all $\ov U\in D(A_1)$, the following equality holds
$$a(U,\ov U)=\suM\iIi \left(\varpsi_1(\li \ov\vx+\ov\u) + \varpsi_2(\li \ov\ux+\ov\v) \right)\ ;$$
by using Lemma \ref{lemma} we obtain $(I-A_1)U=\Psi$ a.e.
.

We can conclude that $A_1$ is m-dissipative, then it is the generator of a contraction semigroup in $Y$, $\mathcal T_1(t)$.

From now on, we follow the path of the proof of Proposition 4.3.3 in \cite{CH} (working in $E$ in place of $X$).
We introduce the set
$$\dsp B_M=\{U\in C([0,T];E):\sup_{t\leq T} \Vert U(t)
\Vert_E
\leq M\}$$
equipped with the distance generated by the norm of $C([0,T];E)$;
 we will find the solution to problem (\ref{A1}) as the unique fixed point in $B_M$ of the function 
$$\dsp\Phi(U)=
\Phi_U(t)=\mathcal T_1(t) U_0+\iot \mathcal T_1(t-s) F(s,U(s))\ ds \ . $$
$\Phi_U\in C([0,T];E)$; moreover, thanks to  the Lipschitz continuity of $F$  in $E$,
for $U\in B_M$, we have
$$\dsp\Vert \Phi_U(t)\Vert_{E}\leq \Vert U_0\Vert_{E}+T L_F(K) M\leq M\ $$
and, for $V\in B_M$,
$$\dsp\Vert \Phi_U(t)- \Phi_V(t)\Vert_{E}\leq L_F(K)\iot \Vert U(t)-V(t)\Vert_{E}
\leq\frac 1 2 \sup_{[0,T]}\Vert U(t)-V(t)\Vert_{E}\ .$$

Then we are able to conclude that $\Phi$ is a contraction in $B_M$ and it has a unique fixed point $U\in B_M$
$$\dsp U(t)=\mathcal T_1(t) U_0+\iot \mathcal T_1(t-s) F(s,U(s))\ ds \ .  
$$ 

Using the above expression we deduce, for $t\in [0,T-h]$, $h>0$,
\be\begin{array}{ll}\dsp
\Vert U(t+h)-U(t)\Vert_Y\leq \mathcal \Vert \mathcal T_1 (h)U_0-U_0\Vert_Y 
+\int_0^h \Vert F(s,U(s)) \Vert_Y ds\\ \\ \dsp
+\iot \left((\Vert f(s)\Vert_{H^1}+\beta) \Vert U(s+h)-U(s)\Vert_Y + \Vert U(s)\Vert_E
\Vert f(s+h)-f(s)\Vert_2\right) ds,
\end{array}\ee
where $\beta=\suM \beta_i$.
Since $f \in C([0,T];H^1(\mathcal A))\cap H^1((0,T);L^2(\mathcal A))$, using Gronwall's lemma, we obtain
$$\Vert U(t+h)-U(t)\Vert_Y\leq C(M,K,T) h\ .$$
Now, using the above inequality and again the assumptions on $f$, we prove that 
$$ \Vert F(s+h,U(s+h))-F(s,U(s)\Vert_2 \leq  \ov C(K,M,T) h.$$
Now we can conclude as in 
Proposition 4.3.9 in \cite{CH}, proving that $U$ is the solution to problem (\ref{A1}), since $U_0\in D(A_1)$.

\end{proof}

\begin{rema}\label{1r}
It is readily seen that the solutions found in the previous two propositions verify
$$\soT\Vert u_t(t)\Vert_2, \soT\Vert v_t(t)\Vert_2,
\soT\Vert \phi_t(t)\Vert_2\leq Q(K,M)\ ,$$
where $Q$ is a quantity depending only on $a,b,\li,\beta_i,D_i$ besides to $M$ and $K$.
\end{rema}

\bt \label{le}{\it (Local existence)}
There  exists a unique local solution $(u,v,\phi)$ to problem (\ref{sysi})-(\ref{compat}),
\be\begin{array}{ll}
(u,v)\in (C([0,T];H^1(\mathcal A))\cap C^1([0,T],L^2(\mathcal A)))^2 ,\\ \\ 
\phi \in C([0,T];H^2(\mathcal A))\cap C^1([0,T],L^2(\mathcal A))
\ .
\end{array}\ee
Moreover, $\phi\in
H^1((0,T);H^1(\mathcal A))$.
\et
\begin{proof}

Let  
$\dsp M>2(\Vert u_0\nhu +\Vert v_0\nhu)$, $K> \Vert \phi_0\nhd+4 M$ 
and $T\leq\min\{ (2 L_F(K))^{-1},1\}$.
We recall that $E:= \dsp \cup_{i\in \M} (H^1(I_i))^2$;
let 

\be B_{MK}=\left\{
\begin{array} {ll}\dsp
(u,v,\phi)\in  (C([0,T];H^1(\mathcal A)))^2\times C([0,T];H^2(\mathcal A)):\\ \\
  \dsp \soT\Vert(u(t),v(t))\Vert_E\leq M\ ,
\soT\Vert \phi(t)\Vert_{H^2}\leq K\ , \\ \\ \dsp
u,\phi\in C^1([0,T];L^2(\mathcal A))\ ,
\soT\Vert u_t(t)\Vert_2, 
\soT\Vert \phi_t(t)\Vert_2\leq Q(K,M) \end{array}
\right\}\ .\ee

We consider the function $G$ defined in $B_{MK}$ :

$$(u^0,v^0,\phi^0)\in B_{MK}\ ,\qquad\ G(u^0,v^0,\phi^0)=(u^{1}, v^{1},\phi^{1})\ ,$$
where $U^{1}=(u^{1}, v^{1})$ is the solution to 
(\ref{A1}) with $f=\phi^0_x$
and $\phi^{1}$ is the solution to problem (\ref{A2}), with $g=a\,u^{1}$.
The previous two propositions ensure that $G$ is well defined from $B_{MK}$ in $B_{MK}$. 
Let 
$$(\hat u^0,\hat v^0,\hat\phi^0)\ ,\ (\ov u^0,\ov v^0, \ov \phi^0)\in B_{MK}\ ,$$
$$(\ov u^{1},\ov v^{1},\ov \phi^{1})=G(\ov u^0,\ov v^0,\ov\phi^0)\ ,\  (\hat u^{1},\hat v^{1},\hat \phi^{1})=G(\hat u^0,\hat v^0,\hat\phi^0)\ ,$$  
$$\ov F=(0,\ov \phi_x^0\ov u^{1}-\beta\ov v^{1})\ ,\ 
\hat F=(0,\hat \phi_x^0\hat u^{1}-\beta \hat v^{1})\ ;$$ 
moreover we denote by $C(M,K)$  constants depending only on the quantities $K,M$ and on the parameters of the problem, and by $\gamma(t)$  functions of $t$ which go to zero when $t$ goes to zero. Then we have

\be\begin{array}{ll}\dsp
\Vert \ov U^{1}(t)-\hat U^{1}(t)\Vert_E
=
 \sup_{[0,T]}\left\Vert\iot \mathcal T_1(t-s)  (\ov F(s)-\hat F(s)) \ ds\right\Vert_E
\\ \\ \dsp
\leq  C(K,M)\ioT\left(\Vert \ov U^{1}(t)-\hat U^{1}(t)\Vert_E
+
\Vert \ov \phi^0(t) -\hat\phi^0(t)\nhd\right)\ dt\  ,
\end{array}\ \ee
whence
\be 
\sup_{[0,T]}\Vert \ov U^{1}(t)-\hat U^{1}(t)\Vert_E \leq \gamma(T)\, C(M,K) \soT \Vert \ov \phi^0(t)-\hat \phi^0(t)\Vert_{H^2}\ ;
\ee
also, using the equations and the above inequality,
\be 
\sup_{[0,T]}\Vert \ov u^{1}_t(t)-\hat u^{1}_t(t)\Vert_2\leq 
C(M,K)\gamma(T)\soT \Vert \ov \phi^0(t)-\hat \phi^0(t)\Vert_{H^2}\ .
\ee

Moreover,  using  (\ref{contr}) and the previous inequalities, we obtain

\be
\dsp\soT
\Vert \ov \phi^{1}(t)-\hat \phi^{1}(t)\Vert_{H^2}
\leq \gamma(T)\,  C(K,M)\sup_{[0,T]}\Vert \ov \phi^{0}(t)-\hat \phi^{0}(t)\Vert_{H^2}
\ee
and finally, using the equation and, again,  the previous inequalities,
$$\dsp \soT \Vert \ov \phi^1_t(t)-\hat \phi^1_t(t)\Vert_{2}\leq \gamma(T)\, C(M,K)\soT \Vert \ov \phi^{0}(t)-\hat \phi^{0}(t)\Vert_{H^2}\ .$$

If  $T$ is sufficiently small, than $G$ is a contraction function in $B_{MK}$ and let $(U,\phi)=(u,v,\phi)$ be  its unique fixed point :
$$ U(t)=\mathcal T_1(t) U_0 +\iot  \mathcal T_1(t-s) F(s,U(s))\, ds\ ,$$
$$ \phi(t)=\mathcal T_2(t) \phi_0 +\iot  \mathcal T_2(t-s) u(s)\, ds\ .$$

Since $u\in C^1([0,T];L^2(\mathcal A))$, arguing as at the end of the proof of Proposition \ref{hA}, we prove that 
$\phi\in H^1((0,T);H^1(\mathcal A))$; thanks to the regularity properties of $\phi$ we can argue as at the end of the proof of Proposition \ref{41} tho prove that $v\in C^1([0,T];L^2(\mathcal A))$. Therefore  $(U,\phi)=(u,v,\phi)$ is the claimed solution.

\end{proof}

\end{section}

\bigskip

\begin{section}
{\bf Global existence }

\bigskip

In this section we prove that, if the initial data are small in a suitable
norm, then  the local solution to problem (\ref{sysi})-(\ref{compat}) given by Theorem \ref{le},
\be\label{spazio}
\begin{array}{ll}
u,v\in \left(C\left([0,T];H^1(\mathcal A)\right)\cap C^1\left([0,T];L^2(\mathcal A)\right)\right)^2\ ,  \\ \\
\phi\in C\left([0,T];H^2(\mathcal A)\right)\cap C^1\left([0,T];L^2(\mathcal A)\right)\cap
H^1((0,T);H^1(\mathcal A))\ ,
\end{array}\ee
can be extended to the time interval $[0,+\infty)$.

We set 
$$\Vert f_i(t)\Vert_2:=\Vert f_i(\cdot,t)\Vert_{L^2(I_i)} ,  \ \ 
\Vert f_i(t)\Vert_{H^s}:=\Vert f_i(\cdot,t)\Vert_{H^s(I_i)}   
.$$

We introduce the functional
\be\label{func}\begin{array}{ll}
\dsp F_T^2(u,v,\phi) := \suM \left(\sup_{t\in[0,T]}  \Vert \u(t)\Vert_{H^1}^2+
\sup_{t\in[0,T]}\Vert \v(t)\Vert_{H^1}^2+
\sup_{t\in[0,T]}\Vert \phii(t)\Vert_{H^2}^2\right) \\ \\
\dsp +\int_0^T 
\left( \Vert u_x(t)\Vert_{2}^2+
\Vert v(t)\Vert_{H^1}^2+\Vert v_t(t)\Vert_{2}^2+\Vert \phi_x(t)\Vert_{H^1}^2
+ \Vert \phi_{xt}(t)\Vert_{2}^2\right)\ dt \ ;\end{array}\ee
we are able to prove that the functional satisfies
the following inequality
\be \label{cubic}
F_T^2(u,v,\phi) \leq c_1 F_0^2(u,v,\phi) +c_2 F_T^3(u,v,\phi) \ , \ c_1,c_2>0\ ,\ee
then,  by standard arguments \cite{Ni}, if $F_0(u,v,\phi)$ is small then $F_T(u,v,\phi)$ remains bounded for all $T>0$.

The  above inequality turns out to be true by some a priori estimates holding for a local solution 
(\ref{spazio}), as we are going to prove in the following propositions. As in the previous section, in order
to simplify the notations, all the proofs are given in the case 
of a graph composed by a single node $N$ and $m$ arcs $I_i$ connecting that node to the external points $a_i$, 
$i\in\M=\{1,2,...,m\}$. 

\bp\label{I}
Let $(u,v,\phi)$ be a local solution (\ref{spazio}) to problem (\ref{sysi})-(\ref{compat}); then

\be\label{uv1}\begin{array}{ll}
\displaystyle
\suM\left( \soT\Vert\u(t)\nld^2+\soT\Vert\v(t)\nld^2 + \beti \ioT \Vert\v(t)\nld^2 dt \right) \\ \\ \leq
C\dsp
\suM\left(\Vert\uo\nld^2+\Vert\vo\nld^2 
+ \soT\Vert\u(t)\nhu \ioT\left(\Vert\phx(t)\nld^2 +\Vert\v(t)\nld^2
\right)dt  \right)\
\end{array}\ee
for  a suitable positive constant $C$.
\epr

\begin{proof}

We multiply the first equation  in (\ref{sysi}) by $\u$, the second one by $\v$ and we sum them; after summing up for $i\in\M$, for $\tau\leq T$ we obtain
\be\begin{array}{ll}
\displaystyle
\suM\left( \Vert\u(\tau)\nld^2+\Vert\v(\tau)\nld^2 + \beti \iota \Vert\v(t)\nld^2 dt \right) \\ \\ \leq
\displaystyle  -C_1 \iota\left(\suE\lai \u(N,t)\v(N,t) -\suU \lai \u(N,t)\v(N,t)\right)dt \\ \\+
C_2\dsp
\suM\left(\Vert\uo\nld^2+\Vert\vo\nld^2 
+ \soT\Vert\u(t)\nhu \ioT\left(\Vert\phx(t)\nld^2 +\Vert\v(t)\nld^2
\right)dt  \right)\
\end{array}\ee
for   suitable positive constants $C_1,C_2$. The transmission conditions (\ref{TC}) imply that the term at node 
is non positive (see Section 3), then the claim.

\end{proof}

\bp\label{II}
Let $(u,v,\phi)$ be a local solution (\ref{spazio}) to problem (\ref{sysi})-(\ref{compat}); then
\be\begin{array}{ll}\displaystyle
\suM\left(\soT
\Vert \vx(t)\nld^2+\soT\Vert\vt(t)\nld^2 +\ioT \Vert\vt(t)\nld^2\, dt\right) 
\\ \\ \displaystyle \leq C \left(\Vert v_0\nhu^2+\Vert u_0\nhu^2\Vert \phi_0\nhd^2\right)
\\ \\ +
 \displaystyle C \suM
\soT\Vert \u(t)\nhu \ioT \left(\Vert\phi_{ixt}(t)\Vert_2^2+\Vert\vt(t)\Vert_2^2\right)\, dt\\ \\
\dsp +C\suM
\soT\Vert \phi_{x}(t)\nhu
\ioT\left(\Vert\vt(t)\Vert_2^2+\Vert\v(t)\nhu^2\right)\, dt
\end{array}\ee
for a suitable positive constant $C$.
\epr
\begin{proof}
In order to avoid the presence of the traces at the node $N$ of the functions $\ut,\vt,\pht$ in the following calculations, we must not derive in time the equations, so we introduce the difference
$\del f(x,t)=f(x,t+h)-f(x,t)$; we have, for $i\in \M$,

\begin{equation}
\left\{ \begin{array}{l}
\left(\del \ut +\lambda_i\del \vx\right)\del\u=0\ ,
\\ \\ \left(\del\vt+\lambda_i \del \ux\right)\del \v= \left(\del(\u  \phx) -\beta_i \del\v\right)\del\v\ .
\end{array}\right.  
\end{equation}
  Summing  the above two equations and  integrating over
$\Ii\times(\delta,\tau)$, for  $0<\delta<\tau< T$,  $\vert h\vert \le \min\{\delta, T-\tau\}$ we obtain

\be\begin{array}{ll}\label{1}\displaystyle
\int_\delta^{\tau} \iIi
\pa_{t}\left(\frac {(\del\u)^2 +(\del \v)^2}2\right) dx \,dt+ 
\int_\delta^{\tau} \iIi\lai\pa_x\left( \del\v\del\u\right)\ dx\, dt \\ \\ \displaystyle
= \int_\delta^{\tau} \iIi\left( \del(\u\phx)\del\v -\beta_i(\del\v)^2\right)\ dx \,dt\ .
\end{array}
\ee

Using the boundary conditions (\ref{vbc}) and the transmission conditions (\ref{TC}) we can compute
$$\dsp\suM
\int_\delta^{\tau} \iIi\lai\pa_x\left( \del\v\del\u\right)\ dx\, dt=
\frac 1 2 \int_\delta^\tau \sum_{i,j\in\M} K_{ij}\left(\del u^j(N,t)-\del\u(N,t)\right)^2 \geq 0\ .$$
Now we divide 
  the equalities (\ref{1}) by $h^2$ , we sum them for $i\in\M$ and,  letting first $h$ and then  $\delta$ go to zero, we obtain
\be\begin{array}{ll}\displaystyle
\suM\left(
\Vert \vx(\tau)\nld^2+\Vert\vt(\tau)\nld^2 +\beta_i \iota\Vert\vt(t)\Vert_2^2\right)\, dt
\\ \\ \displaystyle \leq C_1 
\suM \Vert\vox\nld^2+\Vert\vt(0)\nld^2 
\\ \\ \displaystyle +C_2 \suM
\soT\Vert \u(t)\nhu \ioT\left(\Vert\phi_{ixt}(t)\Vert_2^2+\Vert\vt(t)\Vert_2^2\right)\, dt\\ \\
\dsp +C_3
\soT\Vert \phi_{x}(t)\nhu
\ioT\left(\Vert\vt(t)\Vert_2^2+\Vert\v(t)\nhu^2\right)\, dt
\end{array}\ee
for  suitable positive constant $C_i$.

\end{proof}

\bp\label{III}
Let $(u,v,\phi)$ be a local solution (\ref{spazio}) to problem (\ref{sysi})-(\ref{compat}); then

\be\begin{array}{ll}\displaystyle
\suM\soT\Vert \ux(t)\nld^2\leq C\suM\left(\soT\Vert\vt(t)\nld^2+\soT\Vert\v(t)\nld^2\right)\\ \\
\displaystyle
+C \suM\soT\Vert \u(t)\nhu\left(\soT\Vert\ux(t)\nld^2 +\soT\Vert \phx(t)\nld^2\right)
\end{array}\ee
for a suitable positive constant $C$.
\epr

\begin{proof}

We multiply the second equation by $\ux$, we integrate over $\Ii$ and we sum for
$i\in\M$; using the Cauchy-Schwartz inequality, we obtain the claim
 
\end{proof}

\bp \label{IV}
Let $(u,v,\phi)$ be a local solution (\ref{spazio}) to problem (\ref{sysi})-(\ref{compat}); then

\be\begin{array}{ll}\displaystyle
\suM\ioT\Vert \ux(t)\nld^2\, dt\leq C\suM\ioT\left(\Vert\vt(t)\nld^2+\Vert\v(t)\nld^2\right)
\,dt
\\ \\
\displaystyle
+C \suM\soT\Vert \u(t)\nhu
\ioT\left(\Vert\ux(t)\nld^2 +\Vert \phx(t)\nld^2\right)\,dt
\end{array}\ee
for a suitable positive constant $C$.
\epr

\begin{proof}

We multiply the second equation by $\ux$, we integrate over $\Ii\times(0,T)$ and we sum for
$i\in\M$; using the Cauchy-Schwartz inequality, we obtain the claim.

\end{proof}

\bp\label{V}
Let $(u,v,\phi)$ be a local solution (\ref{spazio}) to problem (\ref{sysi})-(\ref{compat}); then

\be\begin{array}{ll}\displaystyle
 \suM\ioT\Vert \vx(t)\nld^2 \,dt\leq 
C \suM\left(\Vert \vo\nld^2+ \Vert \uo\nhu^2\left(1+\Vert \pho\nhu^2\right)
\right)\\ \\ \dsp  
+C \suM\left(\soT\Vert \vt(t)\Vert_2^2+
\ioT\Vert v_{it}(t)\nld^2\ dt\right)\\ \\
\displaystyle
\displaystyle
+C\suM 
\left(\soT\Vert \u(t)\nhu+\soT\Vert\phx(t)\nhu\right)
\ioT\left(\Vert \v(t)\nhu^2 + \Vert \phi_{ixt}(t)
\nld^2\right)\,dt

\end{array}\ee
for a suitable positive constant $C$.
\epr

\begin{proof}

Using the same notations of the proof of Proposition \ref{II}, by the second equation in (\ref{sysi})
we obtain, for $0<\delta<\tau<T$, $\vert h\vert \le \min\{\delta, T-\tau\}$,
\be\label{2}\begin{array}{ll}\displaystyle
\int_\delta^\tau \iIi\left((\v\del \v)_t-\vt\del \v -\lai \vx\del\u\right)\ dx\ dt\\ \\ \dsp
+\int_\delta^\tau \iIi\lai\left(\v\del\u\right)_x
=
\int_\delta^\tau\iIi \v\left(\del(\u\phx)-\beta_i \del\v\right)\ dx \ dt\ .\end{array}\ee

Using conditions (\ref{vbc}) and (\ref{TC}) 
 we can write
\be\begin{array}{ll}\displaystyle
\lim_{h\to 0}\suM \frac 1 h \int_\delta^\tau \iIi 
\left(-\lai \vx\del\u
+
\beta_i \v\del\v\right)\ dx\ dt  \\ \\ \dsp
=
\lim_{h\to 0}\frac 1 h\suM
\iIi\left( -\v(x,t)\del \v(x,t)\ dx  +\v(x,0)\del \v(x,0)\right)\ dx  \\ \\ \dsp
+\lim_{h\to 0}\frac 1 h\suM\int_\delta^\tau \iIi\left(\vt\del \v 
+ \v\left(\del(\u\phx)\right)\right) dx \ dt\\ \\ \dsp 
-\lim_{h\to 0}\frac 1 h
\sum_{i,j\in\M} \frac {K_{ij}} 2\left(u_j(N,t)-\u(N,t)\right)
\del \left(u_j(N,t)-\u(N,t)\right)\ dt\ .
\end{array}\ee

In order to treat the last term   we set
$$H(t)=u_j(N,t)-\u(N,t)\ $$
and using the continuity of the above function  we have
\be\label{3}\begin{array}{ll}\dsp
\lim_{h\to 0}\frac 1 h \int_\delta^\tau H(t)\del H(t)\ dt
= \lim_{h\to 0}\frac 1 {2h} \int_\delta^\tau 
H(t) \left(\del H(t)-\Delta^{-h} H(t)\right)\ dt \\ \\ \dsp
= \lim_{h\to 0}\frac 1 {2h} \left(-
\int_{\delta-h}^\delta H(t)H(t+h)\ dt + \int_{\tau-h}^{\tau} H(t)H(t+h)\ dt \right) \\ \\ \dsp
=\frac 1 2 \left(H^2(\tau)-H^2(\delta)\right) \ ,\end{array}\ee

Now 
we obtain the claim letting $h$ and then $\delta$ go to zero in (\ref{2}), using equality (\ref{3}).
\end{proof}

\bp\label{VII}
Let $(u,v,\phi)$ be a local solution (\ref{spazio}) to problem (\ref{sysi})-(\ref{compat}); then
\be\begin{array}{ll}\dsp
\suM\left(\soT\Vert \pht(t)\nld^2 +\ioT\left(\Vert \pht(t)\nld^2+\Vert \phtx(t)\nld^2\right)
\,dt\right)\\ \\ \leq
\dsp
C\suM \left(\Vert \pho\nhd^2+\Vert \uo\nld^2+ \ioT\Vert \vx(t)\nld^2\,dt\right)\ ,
\end{array}\ee
for a suitable positive constant $C$.
\epr
\begin{proof}
Using the same notations of the proof of Proposition \ref{II}, by the third equation in (\ref{sysi}) we obtain, for all $i\in \M$,
\be\label{ss}\begin{array}{ll}\dsp
\suM\iIi \int_\delta^\tau \frac {\left( (\Delta^h \phii)^2\right)_t} 2 \ dx\, dt=\\ \\ \dsp
\suM D_i \int_\delta^\tau \iIi\left((\Delta^h \phx) (\Delta^h \phii)\right)_x \, dx\,dt
- \suM D_i \int_\delta^\tau\iIi(\Delta^h\phx)^2
\, dx\,dt\\ \\ \dsp
+ \suM \int_\delta^\tau\iIi\left( a\Delta^h \u) (\Delta^h\phii) -b(\Delta^h\phii)^2\right) \ dx\ dt\ .
\end{array}\end{equation}
We denote by $B^N$ the first term on the right hand side; using the transmission conditions
(\ref{KC}) we have
\be\begin{array}{ll}
\dsp
B^N= \int_\delta^\tau \left(\suE D_i\Delta^h\phii\Delta^h\phx(N,t) - \suU D_i\Delta^h\phii\Delta^h\phx(N,t)\right)\ dt
\\ \\ =\dsp
-\frac 1 2 \sum_{j,i\in\M} \int_\delta^\tau \alpha_{ij} 
\left(\Delta^h\phi_j(N,t)-\Delta^h\phi_i(N,t)\right)^2\, dt\ .
\end{array}\ee
We divide equation (\ref{ss}) by $h^2$ and using 
the Cauchy-Schwartz inequality and letting, first $h$ and then $\delta$,  go to zero  we obtain
\be\begin{array}{ll}\dsp
\suM\Vert \pht(\tau)\nld^2 +\suM\iota\left(\Vert \pht(t)\nld^2+\Vert \phtx(t)\nld^2\right)\,dt \\ \\ \leq
\dsp
C\suM \left(\Vert \pht(0)\nld^2 +\ioT\Vert \ut(t)\nld^2\,dt\right)\, dt\ .
\end{array}\ee
\end{proof}

\bp\label{VIII}
Let $(u,v,\phi)$ be a local solution (\ref{spazio}) to problem (\ref{sysi})-(\ref{compat}); then
\be\begin{array}{ll}\dsp
\suM\left( \soT\Vert \phxx(t)\nld^2 +\soT\Vert \phx(t)\nld^2\right) \leq 
C\suM \left( \soT\Vert \pht(t)\nld^2 +\soT\Vert \u(t)\nld^2\right) \ .
\end{array}\ee
\epr
\begin{proof}

We multiply the third equation in (\ref{sysi}) by $D_i\phxx$, then we sum over $i\in\M$ and using the Cauchy-Schwartz inequality and the boundary conditions we obtain
\be\begin{array}{ll}\dsp
\suM\left( \Vert \phxx(\tau)\nld^2 +\Vert \phx(\tau)\nld^2\right) \leq 
C\suM \left( \Vert \pht(\tau)\nld^2 +\Vert \u(\tau)\nld^2\right) \\ \\ \dsp

+C \left(\suE D_i\phii\phx(N,\tau) - \suU D_i\phii\phx(N,\tau)\right) \, dt\ ;
\end{array}\ee
again, using the transmission conditions (\ref{KC}),
we are able to show that the last term is non positive, so we have the claim.

\end{proof}

\bp \label{VI}
Let $(u,v,\phi)$ be a local solution (\ref{spazio}) to problem (\ref{sysi})-(\ref{compat}); moreover, for all $\nu\in\mathcal P$ let $K_{kj}^\nu\neq 0$ in (\ref{TC}) for, at least, one $k \in \M^\nu$ and for all $j\in\M^\nu$. Then

\be\begin{array}{ll}\displaystyle
\suM\ioT\left(\Vert \phx(t)\nld^2+
\Vert \phxx(t)\nld^2\right)\,dt\\ \\\leq 
\displaystyle
C 
\suM
\ioT\left(\Vert\ux(t)\nld^2+\Vert \v(t)\nhu+\Vert \pht(t)\nld\right) \,dt
\  .
\end{array}\ee
\epr

\medskip

This proposition is crucial in treating the nonlinear terms  $\phx \u$, since it provides an $L^\infty-$ norm's estimate for $\phx$; in the proof, which is given (as for the other propositions) in the case of a single node, we will be leaded to consider
the transmission term
$$
\iota \suE  \left(D_i\u(N,t)\phx(N,t) - \suU  D_i\u(N,t)\phx(N,t) \right)dt
$$
which cannot be discarded by means of sign properties,
 so we are going to use conditions at node to obtain a suitable estimate for it. 
In the case of two arcs, conditions (\ref{TC}) imply that the quantity 
$u_2(N,t) -u_1(N,t)$ is proportional to $v_1(N,t)$ (and to $v_2(N,t)$); then, thanks to conditions (\ref{TC}), we have
\be\label{fxu}\begin{array}{ll}\dsp
\left\vert \ioT( D_1 u_1\phi_{1x}(N,t) -D_2 u_2\phi_{2x}(N,t))\, dt \right\vert
\\ \\ \dsp \leq C \ioT \Vert v_1(t)\Vert_\infty \Vert \phi_{1x}\Vert_\infty\, dt
\leq C 	\ioT\left( \frac {\Vert v_1(t)\Vert_{H^1}^2}{2\ep}+\frac{\ep \Vert \phi_{1x}(T)\Vert_{H^1}^2} 2 \right)\, dt\ .
\end{array}\ee

The further conditions on coefficients $K_{jk}$ assumed in the claim of  the Proposition are necessary to extend this estimate to the cases when the number of arcs is greater than two. We need the following lemma.

\bl
Let $(u,v,\phi)$ be a local solution (\ref{spazio}) to problem (\ref{sysi})-(\ref{compat}) in the case of a single node $N$; moreover,  let $K_{kj}\neq 0$ in (\ref{TC}) for, at least, one $k \in \M$ and for all $j\in\M$. Then, for suitable $\theta_i^j$,
$$\dsp u_j(N,t)=u_k(N,t)+\sum_{i\neq k} \theta_i^j \v(N,t)\qquad \t{ for all } j\in\M.$$
\el
\begin{proof}
Let $\gamma_i=1$ for $i\in \U$ and $\gamma_i=-1$ otherwise; for $i\neq k$ we consider the $m-1$ transmission relations

\be\begin{array}{ll} \dsp
\gamma_i\lambda_i\v(N) =\sum_{j\in \M, j\neq i} K_{ij} (u_j(N)-\u(N)) \\ \\ =\dsp
\sum_{j\in \M, j\neq i,k} K_{ij} (u_j(N)-u_k(N)) -\left(\sum_{j\in \M, j\neq i} K_{ij}\right)
(u_i(N)-u_k(N))
\end{array}\ee
which constitute a linear system in the unknowns $(u_j-u_k)$, $j\neq k$. The assumptions
on $K_{kj}$ ensure that the matrix of the coefficients is non singular (if $k=1$ it is immediate to check that it has strictly dominant diagonal); then the claim follows.
\end{proof}

Let notice that it is never possible to write  $\u(N,t)$ for all $i\in\M$ as linear combination of $\v(N,t)$, since the matrix of the linear system (\ref{TC})  is singular .

\begin{proof} {\it (Proposition \ref{VI})}

We  multiply the third equation in  (\ref{sysi})  by $D_i\phxx$; after summing
for $i\in\M$ , using the Cauchy-Schwartz inequality and the boundary conditions, we obtain

\be\begin{array}{ll}\displaystyle
\suM\ioT\left\Vert \phx(t)\nld^2+\Vert \phxx(t)\nld^2\right)\, dt\leq 
\displaystyle
C 
 \suM \ioT\left(\Vert\ux(t)\nld^2+\Vert\pht(t)\nld^2\right)\,dt
\\ \\
\displaystyle
+
C_1\suM D_i\ioTa\iIi \left(\left(-a\u+b\phii\right)\phx\right)_x\ dx \ dt\\ \\ \dsp
=
C \suM \ioT\left(\Vert\ux(t)\nld^2
+\Vert\pht(t)\nld^2\right)\,dt
\displaystyle
+C_1 \sum_{i,j\in\M}\alpha_{ij}
\ioTa -b \left(\phi_j-\phii\right)^2(N,t)\ dt
\\ \\ \dsp -C_1 a\iota\left( \suE  D_i\u(N,t)\phx(N,t) - \suU  D_i\u(N,t)\phx(N,t) 
\right)\ dt \ ,
\end{array}\ee
where we used the transmission conditions (\ref{KC}).
Thanks to the last lemma we have
\be\begin{array}{ll}\dsp
\suE  D_i\u(N,t)\phx(N,t) - \suU  D_i\u(N,t)\phx(N,t) 
\\ \\ \dsp = \suE  
\left(u_k(N)+ \sum_{j\neq k} \theta^i_j v_j(N)\right)D_i\phx(N,t)\\ \\ \dsp -
 \suU  \left(u_k(N)+ \sum_{j\neq k} \theta^i_j v_j(N)\right)
D_i\phx(N,t) \ ;\end{array}\ee
now we use condition (\ref{FC}) to discard the terms containing $u_k(N)$ and we obtain an estimate similar to (\ref{fxu}). Then the claim follows.

\end{proof}

Now we notice that 
arranging the results obtained in this section, the quadratic terms (not involving the initial data) in the estimates  of Propositions \ref{I} - \ref{VI} 
can be bounded by means  of cubic ones; so, collecting all the estimates, we obtain the inequality (\ref{cubic}) holding for the functional $F$ introduced at the beginning of the section, for all $T>0$. Then we can conclude, as in \cite{Ni,HN}, that, for $F_0$ suitable small , $F$ remains bounded for all $T>0$; this fact proves the following theorem.

\bt (Global existence) Let $K_{kj}\neq 0$ in (\ref{TC}) for (at least) one $k \in \M$ and for all $j\in\M$, $j\neq k$. There exists $\ep_0>$ such that, for 
$$\Vert u_0\Vert_{H^1},\Vert v_0\Vert_{H^1}, \Vert \phi_0\Vert_{H^2}\leq \ep_0 \ ,$$
there exists a unique global solution $(u,v,\phi)$ to problem (\ref{sysi})-(\ref{compat}), 
$$(u,v)\in C([0,+\infty);H^1(\A))\cap C^1([0,+\infty);L^2(\A)) \ , $$
$$\phi\in C([0,+\infty);H^2(\A))\cap C^1([0,+\infty);L^2(\A))
\cap
H^1((0,+\infty);H^1(\mathcal A))
\ .$$
Moreover, $F_T(u,v,\phi)$ is bounded, uniformly in $T$.
\et

\end{section}

\end{document}